\documentstyle{article}
 \def\vt{t\kern-0.22em\raise.18ex\hbox{\char'47}\lower.18ex\hbox{}\kern-0.08em}
 
\newtheorem{th}{Theorem}[section]
\newtheorem{co}{Corollary}[section]

\newtheorem{ob}{Observation}[section]

\newtheorem{rem}{Remark}[section]

\newcommand{\old}[1]{{}} 
\newcounter{obr}
\newcounter{tabul}
\begin{document}
\title{Unit Incomparability Dimension and Clique Cover Width in Graphs\\}
\author{Farhad Shahrokhi\\
Department of Computer Science and Engineering,   UNT\\
farhad@cs.unt.edu
}

\date{}
\maketitle
\thispagestyle{empty}
\date{} \maketitle

%%%%%%%%%%%%%%%%%%%%%%%%%%%%%%%%%%%%%%%%%%%%%%%%%%%%%%%%%%%%%%%%%%%%%%%%%%%%%%%% 

%%%%%%%%%%%%%%%%%%%%%%%%%%%%%%%%%%%%%%%%%%%%%%%%%%%%%%%%%%%%%%%%%%%%%%%%%%%%%%%% 
\begin{abstract}
For a clique cover $C$ in the undirected graph $G$, 
the {\it clique cover graph} of $C$ 
is the  graph obtained by contracting the vertices of  each
clique in $C$ into a single vertex. The {\it clique cover width} of  $G$,
denoted by $CCW(G)$, is the minimum  value of the  bandwidth of 
all clique cover graphs in $G$.  
Any $G$ with $CCW(G)=1$ is 
known to be an incomparability  graph, and hence is called, a {\it unit incomparability graph}.  
We introduced the  {\it unit incomparability dimension of $G$}, denoted by
$Udim(G)$, to be the smallest integer $d$ so that there are 
unit incomparability  graphs  $H_i$ with $V(H_i)=V(G), i=1,2,...,d$, 
so that  $E(G)=\cap_{i=1}^d E(G_i)$. 
We prove a decomposition theorem establishing the inequality  
$Udim(G)\le CCW(G)$.  
Specifically, given any $G$, there are unit incomparability
 graphs  $H_1,H_2,...,H_{CC(W)}$ with $V(H_i)=V(G)$ so that
and $E(G)=\cap_{i=1}^{CCW} E(H_i)$. 
In addition,  $H_i$ is co-bipartite, for 
$i=1,2,...,CCW(G)-1$. Furthermore, we observe that $CCW(G)\ge s(G)/2-1$, where $s(G)$ is the number of leaves in a largest induced star of $G$ , and use Ramsey Theory to give an upper bound on $s(G)$, when $G$ is represented as an intersection graph using our decomposition theorem. Finally, when $G$ is an incomparability  graph we prove that $CCW (G)\le s(G)-1$.

\end{abstract}

\section{Introduction and Summary}
Throughout this paper, $G=(V(G), E(G))$ denotes  a graph
on $n$ vertices. $G$ is assumed to be undirected, unless stated otherwise. 
 The complement of $G$ is denoted by $\bar G$. 
A set $C$ of vertex disjoint   cliques in $G$ is a 
{\it clique cover}, if every vertex in $V$ is, precisely,  
in one clique of $C$. 
Let $L=\{v_0,v_2,...,v_{n-1}\}$ be a linear ordering of vertices
in $G$. The {\it width} of $L$, denoted by $W(L)$, is 
$max_{v_iv_j\in E(G)}|j-i|$. 
The  {\it bandwidth} of $G$ \cite{CC}, \cite{DP}, 
denoted by $BW(G)$,  is the smallest  width of all linear orderings
of $V(G)$. The  bandwidth problem  is   well studied
and  has intimate connections to other important  concepts
including graph separation \cite{Bo}, \cite{BK}, \cite{Hal}. 
Unfortunately, computing
the bandwidth is $NP-hard$, even if $G$ is a tree \cite{GJ}. 
\vskip .4cm

In \cite{Sh2} we introduced the {\it clique cover width problem}
which is a  generalization of the bandwidth problem.
For $C$ a clique cover of $G$, 
let the {\it clique cover graph}  of $C$, denoted by $G(C)$,  be the graph 
obtained from $C$ by contracting the vertices in each clique to one vertex in $G$. 
Thus $V(G(C))=C$, and $E(G(C))=\{XY| X,Y\in C, there~is~xy\in E(G)~with~x\in X,~and, y\in Y\}$.
The {\it clique cover width} of  $G$, denoted by $CCW(G)$, 
is the minimum value of $BW(G(C))$,
where the minimum is taken over all clique covers  $C$ of $G$.
Note that  $CCW(G)\le BW(G)$, since $\{\{x\}|x\in V(G)\}$ is a  
trivial clique cover in $G$  containing  cliques  that have only one  
vertex. We highly suspect that the problem of computing
CCW(G) is N P − hard, due to the connection to the bandwidth prob-
lem.
Let $C$ be a clique cover in $G$. Throughout this paper,   we will  write 
$C=\{C_0,C_1,...,C_t\}$ to indicate that $C$ is an ordered set of cliques.
For a clique cover $C=\{C_0,C_1,...,C_t\}$, in $G$, let
the {\it width}  $C$, denoted by $W(C)$, denote
$\max\{|j-i|~|C_iC_j\in E(C(G))\}$.  Observe that,
$W(C)= \max\{|j-i||xy\in E(G), x\in C_i,y\in C_j,C_i,C_j\in C\}$. 
\vskip .4cm

An crucial tool for the design of a divide and
conquer algorithm is separation.
The planar separation  \cite{LT}
asserts  that any $n$ vertex planar graph can be separated into two
subgraphs, each having at most $2n/3$ vertices,   by removing $O({\sqrt n})$ vertices. 
The key application of the clique cover width 
is in the derivation of separation theorems in graphs, where separation can be
defined for  other types {\it measures }\cite{Sh1}, instead of just the number of vertices. 
For instance, given a  clique cover $C$ in $G$,
can $G$ be separated  by removing a *small* number of cliques in  $C$ so that
each the two remaining subgraph of $G$ can be covered by at most ${\alpha |C|}$ cliques from $C$, 
where $\alpha<1$ is a constant \cite{Sh1,Sh3}?      
\vskip .4cm

Any (strict) partially  ordered set  $(S,<)$ has a directed acyclic  graph 
$\hat G$ associated with it in a natural way: $V(G)=S$, and $ab\in E(G)$
if and only if $a<b$. The {\it comparability graph}  associated with $(S,<)$
is  the undirected graph which is obtained by dropping the orientation on
edges of $\hat G$ \cite{C2,To2}.  
The complement of a comparability graph is an  {\it incomparability graph}. 

Any graph G with $CCW (G) = 1$ is known to be an incomparability graph \cite{Sh2}, and hence we call such a G a unit incomparability graph.
Clearly, any co-bipartite graph, or the complement of a bipartite graph is
a unit comparability graph. In addition, it is easy to verify that any unit
interval graph is also a unit incomparability graph. Thus, the class of unit
incomparability graphs is relatively large.
%It is well known that incomparability graph have geometric representation 
%in the plane \cite{C1}. In addition,  some of earliest significant 
%advances in  geometric graph  theory was achieved by representing
%any geometric graph as the intersection graph of a four  incomparability 
%graphs \cite{PT}.
%Moreover,  recent  studies on the structure of planar curves 
%has identified large incomparability graphs that appear 
%as subgraphs of  intersection graph of strings \cite{FP2}.
%These result highly suggest that incomparability graphs may be viewed as 
%the  "building blocks" of many geometric structures. 

Let $d\ge 1$, and for $i=1,2,...,d$
let $H_i$ be a graph with $V(H_i)=V$, and  let $G$ be a graph with
$V(G)=V$ and 
$E(G)=\cap_{i=1}^dE(G_i)$. Then,  we say that  $G$ is the 
{\it intersection graph}
of $H_1,H_2,...,H_d$, and write $G=\cap_{i=1}^dH_i$ \cite{Su1}.  
Let the incomparability dimension of G, denoted by $Idim(G)$, be the smallest integer d so that there are
d incomparability graphs whose intersection graph is G. Similarly, let the
unit incomparability dimension of G, denoted by $Udim(G)$, be the smallest
integer d so that there are d unit incomparability graphs whose intersection
is G. In this paper we focus on the connection between the the clique cover
width and the unit incomparability dimension of a graph. Our work gives
rise to a new way of representing any graph as the intersection graph of
unit incomparability graphs.

Recall that Boxicity, and Cubicity of a graph, denoted $Box(G)$ and $Cub(G)$
respectively, are  the smallest integer   $d$ so that $G$ is 
the intersection graph of $d$ interval graphs, or unit interval graphs, 
respectively \cite{Ro}. Recent work \cite{Su1,Su2}  has elevated the 
importance 
of the Boxicity and the Cubicity, by improving the upper bounds, and  by linking these concepts to 
other important graph parameters including  treewidth 
\cite{Bo}, Poset dimension \cite{To2}, and crossing
numbers \cite{Sh4}.  Clearly, $Idim(G)\le Udim(G)$, $Box(G)\le 
Cub(G)$, $Idim(G)\le  Box(G)$, and $Udim(G)\le Cub(G)$. 
 
While Cubicity
and Boxicity are related to unit incomparability dimension, and while the
results in \cite{Su1,Su2} are extremely valuable, these results do not imply, or
even address the concepts and results presented here, particularly, due to
the focus of our work on the pivoting concept of the clique cover width.
\vskip .4cm

%One may view concept of Boxicity and
%Cubicity are as  models for  representation graphs as the intersection graph
%of incomparability graphs, since  interval graphs are incomparability graphs
%associated with interval orders \cite{To1}. 
%Recent work \cite{Su1,Su2}  has elevated the 
%importance 
%of the Boxicity and the Cubicity, by improving the upper bounds, and  by linking these concepts to 
%other important graph parameters including  treewidth 
%\cite{Bo}, Poset dimension \cite{To2}, and crossing
%numbers \cite{Sh4}. 
%\vskip .4cm

\vskip .4cm
In Section two we prove a decomposition theorem that establishes the
inequality $Udim(G) \le CCW (G)$, for any graph $G$. Furthermore, we observe that $CCW (G)\ge (s(G)/2)-1$, where $s(G)$ is the largest number of
leaves in an induced star in $G$, and use Ramsey Theory to give an upper
bound on $s(G)$, when $G$ is represented as an intersection graph using our
decomposition theorem.
 In Section three we study the clique cover width
problem in incomparability graphs, and prove that $s(G)-1\ge CCW(G)$,
when $G$ is an incomparability graph. The results give rise to polynomial
time algorithms for the construction of the appropriate structures.

{\bf Remark.} 
In work in progress, the author  has   improved the  upper bound on $Udim(G)$ (in the decomposition theorem) to
$O(log(CCW(G))$. This drastically improves all upper bounds presented here.
 
%The result  suggested to explore

%the  connection between the clique cover width,
%and incomparability graphs.

\section{Main Results}   
Now we prove the decomposition result.

\begin{th}\label{t3}(Decomposition Theorem)
{\sl
Let $C=\{C_0,C_1,...,C_t\}$ be a 
a clique cover in $G$.
Then, there are $W(C)$ unit incomparability
graphs $H_1,H_2,...,H_{W(C)}$ 
whose intersection is $G$.
Specifically, 
$H_i$ is a co-bipartite graph, 
for  $i=1,2...,W(C)-1$. Moreover, the constructions
can be done in polynomial time.

}
\end{th}

{\bf Proof.} 
For $i=1,2,...,W(C)-1$, we define
 a graph $\bar H_i$ on the vertex set $V(G)$ and the edge set
$$E({\bar H_i})=\{xy\in E({\bar G})| x\in C_l,y\in C_k,|k-l|=i\},$$
and prove that  
$\bar H_i$ is a bipartite graph.  
For $i=1,2,...,W(C)-1$,
let $odd(i)$ and $even(i)$ denote the set of all
integers $0\le j\le t$, so that $j=a.i+r,r\le i-1$, where $a\ge 0$ is
odd, or even, respectively, and note that $odd(i)\cup even(i)$
is a partition of $\{0,1,2,...,t\}$.
Next, for $i=1,2,...,W(C)-1$, let 
$V_1=\{x\in C_j| j\in odd(i)\}$ and 
$V_2=\{x\in C_j| j\in even(i)\}$, and note $V_1\cup V_2$ is a partition
of $V(G)$ so  that 
any edge in $E({\bar H}_i)$ has one end point in $V_1$ and the other
end point  in $\in V_2$. 
Therefore, for $i=1,2,...,
W(C)-1$,  
$\bar H_i$ is bipartite, and  consequently, $H_i$ is co-bipartite.
Next, let 
$\bar H_{W(C)}$ 
be a graph on the vertex set
$V(G)$ and the edge set 
$$E({\bar H}_{W(C)})=\{xy\in E({\bar G})|
x\in C_l,y\in C_k, |l-k|\ge 
{W(C)}\}.$$
Let $xy\in E({\bar H})$, then
$y\in C_l$  and $x\in C_k$ with  
$|l-k|\ge 
W(C).$
Now orient $xy$ from $x$ to $y$, if 
$l\ge  k+
W(C)$, otherwise, orient $xy$ from $y$ to $x$.
It is easy to verify that this orientation is transitive, and hence
${\bar H}_{W(C)}$ is a comparability graph. 
Consequently, 
$H_{W(C)}$ is an incomparability graph.
We need to show $CCW(H_{W(C)})=1$.  
First, observe that any  consecutive subset  cliques of  in $C$ 
of cardinality at at most $W(C)$  is a clique
in $H_{W(C)}$. Next, let $t=a.W(C)+r, r\le W(C)-1$. For $i=0, 1,...,a-1$,  
let $S_i$ be the set of all $W(C)$ consecutive cliques in $C$ starting at 
$i.W(C)+1$ and ending at $(i+1)W(C)$. Define $S_a$ to be the set of $r$
consecutive cliques in $C$, starting at $C_{a.W(C}+1$ and ending in $C_t$.
It is easy to verify that  $S=\{S_0,S_2,...,S_a\}$ is a clique cover 
in $H_{W(C)}$ with $W(S)=1$, and thus $CCW(H_{W(C)})=1$. 
$\Box$

\begin{rem}\label{r1}
{\sl If $G$ is a clique, then $CCW (G) = 0$, where, $Udim(G)=
1$. In addition, if $G$ is disconnected, then $CCW (G)$ and $Udim(G)$ equal
to the maximum values of these parameters, respectively, taken over all
components of $G$.}
\end{rem}

A Simple consequence of Theorem \ref{t3}  is the following.

\begin{co}\label{c1}
{\sl $CCW (G)\ge Udim(G)$, for any connected graph $G$ which
is not a clique.
}
\end{co}

In light of Theorem \ref{t3} , one may want to have estimates for $CCW (G)$.
Let $s(G)$ denote the number of leaves for a largest induced star in $G$. When
$|V (G)| ≤ 2$, we define $s(G) =1$.

\begin{ob}\label{ob1}
{\sl Let $C=\{C_0, C_1...,C_k\}$ be a clique cover in $G$, then, 
$W(C)\ge \lceil{s(G)\over 2}\rceil-1$.}

\end{ob}
{\bf Proof.} 
Let  $S$ be an induced  star with center $r$ and $s(G)$ leaves. 
Then,  $r\in C_i$ for some $0\le i\le t$.
Note that no two leaves in $S$ can be in the same clique of  $C$, 
and hence there must be an edge $rx\in E(S)$ with  $x\in C_j$ for
some $0\le j\le t$ so that  $|j-i|\ge 
\lceil{s(G)\over 2}\rceil-1$.
Consequently, $W(C)\ge \lceil{s(G)\over 2}\rceil-1$.
$\Box$

Ideally, one would like to have a converse for Theorem \ref{t3}, where the
clique cover width of $G$ can be estimated using the clique cover widths for
the factor graphs. Unfortunately, we have not been to derive a general
result that can used effectively. In addition, Observation \ref{ob1} poses the
problem of finding an upper bound for $CCW (G)$ in terms of $s(G)$, only.
Unfortunately this is also impossible. For instance, take the $n\times n$ planar
gird $G$, then $s(G) = 4$, but $CCW (G)$ is unbounded. (We omit the details.)
Nonetheless, one may pose the related question of finding an upper bound
for $s(G)$ in $G=\cap_{i=1} H_i$ , using $s(Hi )$, $i = 1, 2, ..., d$.

For integers, $n_1,n_2...,n_c$ let $R(n_1,n_2,...,n_c)$ denote the general Ramsey number.  Thus for any 
$t\ge R(n_1,n_2,..., n_c)$, if the edges   of $K_t$ are
 colored with $c$  different colors, then for some $1\le i\le c$,  we always 
get a complete subgraph of $K_t$,   on $n_i$ vertices,  whose edges are 
all colored with color $i$. 

\begin{th}\label{t4}
{\sl 
 Let $G=\cap_{i=1}^d H_i$,  then $s(G) < R(s(H_1 ) + 1, s(H_2 ) + 1, ..., s(H_d ) + 1)$. }
 \end{th}

 {\bf Proof.} Assume to the contrary that $s(G)\ge  R(s(H_1 ) + 1, s(H_2 ) +
1, ..., s(H_d ) + 1)$, and let $S$ be an induced star in G rooted at $r$ with
$s(G)$ leaves. Note that $S=\cap_{i=1}^d S_i$ , where for $i = 1, 2, ..., d$, $S_i$ is an
induced subgraph of $H_i$ with $V (S_i )=V (S)$. Now let $H$ be a complete
graph on the vertex set $V(H) = V (S) − \{r\}$. Thus, the vertex set of
$H$ is precisely the set of leaves in $S$. Now for $i = 1, 2, ..., d$, and for
any $ab\in  E(H)$, assign color $i$ to $ab$, if and only if, $ab \notin  E(S_i )$. 
Since$|V (H)| = s(G)\ge R(s(H_1 ), s(H_2 ), ..., s(H_d ))$, there must be a monochromatic 
complete subgraph of H on color $i$, with least $s(H_i ) + 1$ vertices, for
some $1\le i\le d$. Let $W_i$ denote the set of vertices for this complete sub-
graph, and note that our coloring policy implies that $W_i$ is an independent
set of vertices in $H_i$ . Thus $W_ i\cup \{r\}$ is an induced star in $H_i$ with at least
$s(H_i ) + 1 > s(H_i )$ leaves which is a contradiction. $\Box$

\begin{co}\label{c3}
{\sl  For any $G$, 
$s(G)<R(\underbrace{3, 3, ..., 3}_{CCW(G)-1~times}, 4)$.}
\end{co}

{\bf Proof.} By Theorem \ref{t3} , there are $CCW (G)$ unit incomparability graphs
$H_1 , H_2 , ..., H_{W (C)}$  whose intersection is $G$. Observe, that $s(H_i)\le  2$, since
$H_i$ is co-bipartite for $i = 1, 2, ..., CCW(G)-1$, and, $s(H_{CCW (G)})\le 3$, since
$H_{CCW (G)}$ is a unit incomparability graph. Now apply Theorem \ref{t4}. $\Box$

\section{Clique Cover Width in Incomparability Graphs }

\begin{th}\label{t5}
{\sl Let $G$ be an incomparability  graph, then there is a clique cover
$C=\{C_0, C_1,...,C_k\}$ so that
 $s(G)-1\ge W(C)\ge 
\lceil{s(G)\over 2}\rceil-1$. Moreover, if graph  $\hat G$, or the transitive
orientation of $\bar G$,  is given in its adjacency list
form, then $C$ and $W(C)$ can be computed in $O(|V(G)|+|E(G)|)$ time.
}
\end{th}

{\bf Proof.} 
Let $C=\{C_0, C_1,...,C_k\}$ be a greedy 
clique cover in $G$, where $k+1$ is the size of largest  independent set in $G$.
Thus, $C_1$ is the set of all sources  
in the $\hat G$, 
$C_2$ is the set of all forces that are obtained  after removal of $C_1$,
etc. 
The lower bound for $W(C)$ follow   from Observation \ref{ob1}.   
For the upper bound, let  $e=ab\in E(G)$ with $a\in C_i,b\in C_j, j>i$ so that 
$W(C)=|j-i|$. Now let  $x_j=b$,  then for $t=i,i+1,...,j-1$ there
is $x_t\in C_t$ so that  $x_tx_{t+1}\in E(\hat G)$. It follows that
for $t,p=i,i+1,...,j, p>t$, we have $tp\in E({\hat G})$ and thus 
$tp\notin E(G)$. We conclude that the induced graph on $a,x_i,x_{i+1},...,x_j$
is a star in $G$, with center $a$, having $j-i+1=W(C)+1$ leaves. Consequently,
$s(G)\ge j-i+1\ge W(C)+1$. To finish the proof, and for the algorithm one
can apply topological ordering to  graph $\hat G$ to compute $C$ and $W(C)$
in linear time.
$\Box$

\begin{co}\label{c5}
{\sl  $G$ be an incomparability graph then, $CCW (G)$ can be
approximated within a factor of two, in $O(|V | + |E|)$ time.}
\end{co}

\begin{co}\label{c6}
{\sl  Let $G$ be a incomparability graph, then $Udim(G)\le s(G)-1$.}
\end{co}

\begin{rem}\label{r4}
{\sl  Let $G$ be a co-bipartite which is not a clique. Then, $S(G) = 2$
and $CCW (G) = 1$, and hence the upper bound in Theorem 3.1 is tight. Let
G be graph consisting of 3 cliques $C_1$ , $C_2$ , $C_3$ , in that order, each on 4 vertices, with $x_i\in  C_i, i = 1, 2, 3$ and additional edges, $x_1 x_2$ and $x_2 x_3$ . Then,
$s(G) = 3$, and $CCW (G) =1$, and hence the lower bound in Observation
2.1 is tight.}
\end{rem}

%The author has succeeded to  improve Theorem \ref{t3} to the following.

%\begin{th}\label{t4}(Decomposition Theorem \cite{Sh3})
%{\sl
%Let $C=\{C_0,C_1,...,C_t\}$ be a
%a clique cover in $G$, then there is an integer
%${\hat k}\le 2\lceil\sqrt{W(C)}\rceil$,
%co-bipartite graphs
%$H_i, i=1,2,...,{\hat k}-1$, and
%a unit incomparability  graph $H_k$
%so that $G=\cap_{i=1}^{\hat k}H_i$.
%}
%\end{th}

\end{document}